# Flipping Linear Algebra using a MOOC platform

Ana Moura Santos, *Department of Mathematics*
*Instituto Superior Técnico, University of Lisbon, Portugal*
*ana.moura.santos@tecnico.ulisboa.pt*

Luis Costa, *Programa de Engenharia de Sistemas e Computação,*
*Coppe/UFRJ. Federal University of Rio de Janeiro, Brazil*
*luisfcosta@cos.ufrj.br*

*Abstract*— Between 2017 and 2019, a standard Linear Algebra course from Instituto Superior Técnico, University of Lisbon, used virtual learning content, mainly videos and formative assessment, delivered at the institution's MOOC platform to support a flipped classroom strategy. This strategy has involved around 100 first-year Computer Science enrolled students each year and was implemented by a faculty member of the Department of Mathematics. The flipped classroom was pedagogically developed and underwent new evaluations each year, with positive impacts on the way students and teachers work during the semester, within the context of teaching/learning mathematics. We will present and discuss the improvements and the results of this three-year hybrid experiment, analysing data from students' responses to questionnaires, online course completion rates, and data collected from students' final grades in the Linear Algebra course. Based on the findings, we have reason to believe that similar flipped classroom practices are a good response to the necessary changes in higher education towards more student-centred practices, especially in the areas of STEM higher education.

## 1. Introduction

Information and communication technologies have provided transformations in several areas in society, especially in recent years. Due to these changes, the education area particularly and globally has provided new relationships between teachers and students (Bates, 2015; Hamdan et al., 2013). These new relationships, together with technological tools, such as videos, interactive quizzes, and constant activities, outside the classroom, have placed the student in the central role of learning. One of these forms of learning is known as the flipped classroom (Bishop & Verleger, 2013), which in the context of mathematics has already been analysed extensively, mostly within the secondary education framework, (Cheng-Yu et al., 2019; Zengin, 2017), but as far as we know not in the Portuguese secondary curriculum (Santiago et al., 2012). In this article, we seek to describe the flipped classroom implementation, within a standard undergraduate Linear Algebra (LA) course, making use of virtual learning content structured as an online course made available by the MOOC Técnico platform[1] of

---
[1] https://courses.mooc.tecnico.ulisboa.pt

Instituto Superior Técnico (IST)[2], University of Lisbon. The implementations have run during the fall semesters of three consecutive years from 2017 until 2019. The average attendance was 100 enrolled students per semester. The results obtained so far demonstrate that this type of learning methodology should continue to be encouraged, not only because students appreciate it a lot, but because it allows them to be more engaged with the course content and perform well at the summative assessments.

*1.1. MOOCs delivered in a Virtual Learning Environment*

A Massive Open Online Course (MOOC) is a carefully designed course on a certain subject that can be accessed online and made available to an unlimited number of participants. In this case we are considering here, the use of an institutional MOOC platform that uses open access, together with convenience and scalability makes these online courses an extremely interesting resource for pedagogical experiments. The term MOOC was coined in 2008 by Dave Cormier[3] and has since then been widely used, both for cMOOCs, more organic courses focused around learner-generated content, and xMOOCs which have educator generated content (often in the form of videos). In the last years, MOOCs have experienced a considerable evolution and major education providers offer xMOOCs, that have increasingly incorporated some features from cMOOCs: discussion forums, peer reviews, and tutor's feedback.

The acronym MOOC[4] stands for the following distinctive characteristics, all of them in use in the MOOC Técnico platform:

- Massive: easy scalability since all content is available online and server capacity is being taken into account;

- Open: with direct access and minimal administrative requirements for registration; Open content upon free registration, and free honour certificates;

- Online: through the internet that allows these courses to be largely accessed, also enables the use of several multimedia and interactive resources, the possibility of automated grading, via assessment activities (quizzes) and participants'-based peer reviews;

- Course: the delivered material contains relevant information in an organizsed and structured way. The pace, duration, intensity, role of instructor and tutors, type of used resources and assessment activities are clearly stated in advance.

The online courses launched so far in the institutional MOOC platform, a customised Open edX platform that serves not only as a Virtual Learning Environment (VLE) for the IST students but offers open and free access to everyone, are developed through multidisciplinary collaboration between our teaching staff community and the MOOC graphic and motion design team. All online courses are in Portuguese and/or English and include video, transcripts and texts from other languages. Typically, these MOOCs have a planned duration of 4 to 5 weeks corresponding to a workload of 5 to 8 hours per week (circa 1.5 ECTS per MOOC) since several of them are complementary activities to face-to-face courses, usingblended learning strategies (Gomes et al., 2018). The multimedia solutions for content are mainly short and medium-length educational videos, visually stimulating, and each course foster different forms of assessment: Quizzes, Peer review, and Formative assessment, which allow any enrolled participant to earn a free certificate upon successfully achieving at least 60% of the planed graded activities.

---

[2] Also known as Técnico Lisboa, the largest Portuguese school of Architecture, Engineering, Science, and Technology, integrated in University of Lisbon, with approximately 11,500 students, distributed among 18 different undergraduate degrees, and several MSc. and Ph.D. programs.

[3] The term MOOC was coined in 2008 by Dave Cormier to an open online course designed and led by Georg Siemens (Athabasca University) and Stephen Downes (The National Research Council, Canada).

[4] In this paper we are considering only xMOOCs, using the shorter designation MOOCs.

## *1.2. Flipped classroom strategies*

The common definition of a flipped classroom, according to several authors (Bishop & Verleger, 2013; Brame, 2013; Selingo, 2016), is a collective or individual learning experience based on problem-solving which includes access to resources and materials outside the classroom. One of the most extensively used technological resources, according to these authors, is the use of videos. In this type of learning, the teacher is no longer the centre of learning (Ebert & Culyer, 2013; Hutchings & Quinney, 2015), sharing space with these technological resources and sometimes with the knowledge produced amongst the students (Prud'Homme-Généreux, 2016).

The flipped classroom approach has been used in several Higher Education Institutions (HEI) within the framework of different subjects (Herreid & Schiller, 2013; Moffett, & Mill, 2014; Raffaghelli et al., 2018). As an example of the use of this practice, there are reports of extensive use at Universities such as Michigan State and Case Western Reserve, generating changes in their behaviours and altering learning results (Selingo, 2016).

More recently, as the report of Gwo-Jen et al. ( 2019) highlights, it can be seen that flipped methodologies have moved into a new era, when researchers are more focused on investigating whether adopting innovative strategies or technologies that allow flipped methods to act more efficiently. This is in contrast to earlier studies, when researchers were more stressed about the effectiveness of flipped strategies when compared with performances from more traditional learning methodologies, e.g. in (Berge, 1995; Brame, 2013).

Here we adopted Bishop and Verleger's (2013) definition of the flipped classroom methodology as a technology-supported pedagogy that consists of two components: a) direct individual computer-based instruction outside the classroom through video lectures, and b) interactive group learning activities inside the classroom (see Fig. 1).

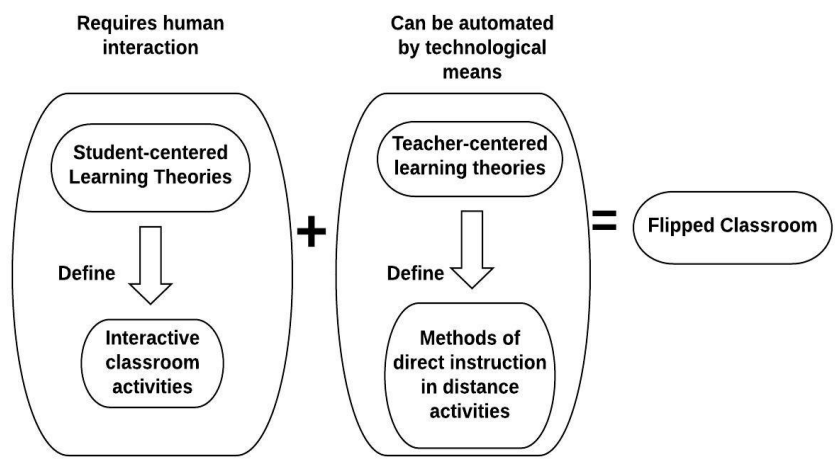

Fig. 1.   Scheme adapted from Bishop and Verleger (2013).

Inverting the classroom, therefore, means for us that students gain initial exposure to new material outside the class, usually through lecture videos or reading, and use class time on face-to-face more complex activities according to Bloom's taxonomy (Bloom, 1956): applying, analysing and evaluating. When coming to campus, students perform the most demanding activities using previously acquired knowledge to solve problems, deduce new relations between concepts, debate and validate original hypotheses (see Fig. 2). In our opinion the study conducted by Ofugi (2016) contributes to a better understanding of the concept of the flipped classroom, identifying the main characteristics of this methodology in opposition to the main features of traditional classes (see Fig. 3).

|  | **Tradicional Classroom** | **Flipped classroom** |
|---|---|---|
| **Teacher preparation time for classes** | Planning the content and preparation of activities to be developed in the classroom | Preparation of asynchronous resources (videos, exercises, simulations) and activities to be developed in the classroom |
| **Student role** | Hear explanations from the teacher in class and then perform consolidation activities at home | Study previously at home and perform activities in the classroom, preferably in groups |
| **Main role of the teacher during face-to-face class** | Content transmitter | Mediator and facilitator |
| **Classroom activities** | Explanation of contents, clarification of doubts and consolidation by solving exercises | Deepening of the content seen at home, clarifying doubts and varied consolidation activities |
| **Activities performed at home** | Various activities to stregthen knowledge on the content seen/heard during the class | Access to asynchronous resources and simple activities to understand the content |
| **Content access location** | First time in the classroom | First through asynchronous resources, which can be accessed in different locations |

Fig. 2.   Table adapted from Ofugi (2016).

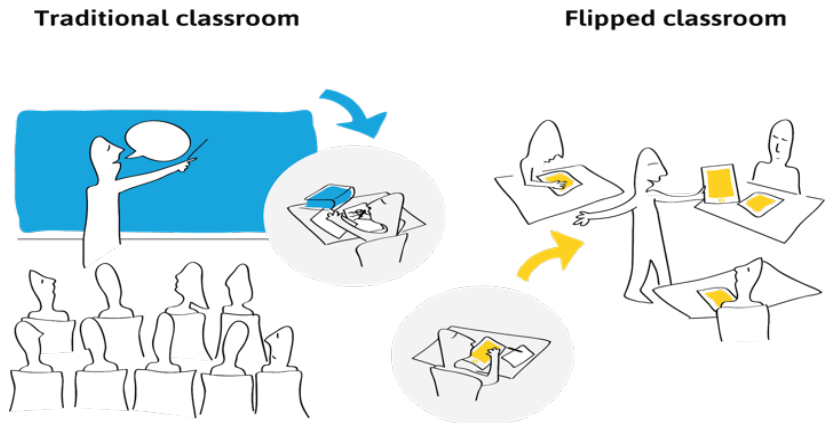

Fig. 3.  Picturing the main differences between a traditional expository lecture and a flipped classroom environment.

Here we aim to analyse firstly in more depth some relevant pedagogical issues about the implementation of the three years flipped classroom methodology, and then discuss the improvements and results of those practices.

**2. Proposal for flipping a Linear Algebra course**

During the fall semester of 2017/2018 academic year, it was presented a proposal for flipping the class within a LA course to the Department of Mathematics of IST. The proposal, supported by the Pedagogic Council, was approved under certain restrictions imposed by the Department of Mathematics (DM), and applied for the first time in that semester to circa 100 enrolled undergraduates (first-year students) in a Computer Science (CS) Engineering Degree. In the subsequent fall semesters of 2018/2019 and 2019/2020 the flipped classroom methodology was again applied to more than 90 first year on-campus students that entered the University in the corresponding academic years. Observe that the LA course is a regular and mandatory course offered in every academic term by the DM to about 1,500 students from different Engineering degrees in IST and therefore if there's more support it can be applied to more students.

*2.1. Flipping Linear Algebra using a MOOC*

At the moment, the LA syllabus includes mastering matrix calculus and methods in solving systems of linear equations, knowledge about vector spaces and linear transformations, which includes the study of canonical forms of matrices, eigenvectors and eigenvalues[5]. In our case, we adopted the textbook "Linear Algebra and its applications" by David Lay (Lay et al., 2016), where the standard LA subjects are supported by applications that can range from Computer Graphics, linear dynamical systems, Markov chains, PageRank algorithm, image compression that is based on SVD, just to name some relevant applications to CS students.

The proposed flipped classroom methodology is mainly based on the content of a MOOC about eigenvalues and eigenvectors, which runs openly in MOOC Técnico platform, during the last five weeks, out of the 14 weeks, of the LA semester. First-year CS students enrolled in LA are all advised to enrol in the MOOC and thereafter are exposed to the course curriculum outside the classroom (watching videos and interacting with MOOC's content and activities). Later that week they shall come to the campus and interact with each other and with the teacher in the classroom. The MOOC is an online course specially designed to fit the last one fourth LA standard course syllabus and is four-five weeks long (6 to 7 hours of estimated effort per week). Besides, the content is also aimed at a broader audience with an interest in current applications of Mathematics and serves well the remaining enrolled participants.

---

[5] The singular value decomposition, well known as SVD is a non-mandatory subject although its importance for STEM students.

The choice for the topic of eigenvalues and eigenvectors, among the LA syllabus, is based on two main reasons. The first one relies on the relevance of this LA subject to engineering applications, many of them already mentioned: Stability of linear dynamical systems, Markov chains, PageRank algorithm, Image compression, Principal Analysis Component (PCA). The second reason is connected to pedagogic reasoning that for first-year students it's better to introduce such innovative learning strategies after they established mutual trust with the teacher during the first weeks of the semester.

Specific questions and queries that arise when trying to submit answers to the online quizzes and assessment activities, should be discussed both in the MOOC forum and during face-to-face classes. Particularly, there was an expectation for more vivid face-to-face classes, where most CS undergraduates engage in problem-solving, or in figuring out connections between concepts, using what they've been learning out-of-classes, rather than spending time in being told "here's a definition", "here's a theorem" (see Table 1).

## *2.2. Challenges of student-centered practices*

We acknowledge that undergraduate students at IST, and in general students issued from Portuguese high-schools, are used to passive classroom experiences with a strong emphasis in summative assessment and do not have much experience on autonomous learning. In fact, the report of Santiago, et al. (2012) highlights how in Portuguese high-schools it's paid little attention to the development of students' capacity to regulate their learning through self- and peer-assessment. In this report, one can find the recommendation that evaluation frameworks in Portugal, in particular within secondary education, the education institutions should focus much more directly on the quality of learning and teaching and their relationships to student outcomes.

As a consequence, the flipped classroom practice proposed here, is at the same time a novelty and a challenge for first-year university students. To prepare them for the late semester challenge, we always feel the need for the first days to foster students' active learning behaviours through out-of-class preparation, participation in group dynamics, and critical thinking, just to mention a few activities. The tutorial/exercises classes with group discussions and comments from peers and tutor on the students writing on the wall (there will be a more detailed description in the next section), together with a formative online assessment (Vonderwell & Boboc, 2013; Moura Santos & Ribeiro, 2017) operating in the first days of the semester, aiming to build a self-confident, more autonomous and related group of students.

## *2.3. A MOOC on a Linear Algebra subject*

As previously described, the MOOC Técnico platform is a VLE where are embedded videos with lecture contents, carefully planned and graphically appealing, and where the content can be released in sequential order with one weekly topic followed by next weekly topic. That is, we took under consideration the running of a structured online course, a cMOOC, organized around an HEI standard curriculum, with content mainly delivered through videos (Raffaghelli et al., 2018). The other pedagogical elements of the online course, mainly quizzes and automated testing (assessment activities), induce a "linear journey" for the learners and foster moments for online discussion (course forum).

On top of these activities, one can also find in the MOOC on eigenvalues and eigenvectors interactive tools, interactive demonstrations, for practicing concepts and applications, and a glossary summarizing key concepts. This MOOC design works very well for a fundamental mathematical syllabus as LA, where one builds the present version based on the knowledge from previous "lectures" and serves an extended group of STEM stakeholders.

On the other hand, the described virtual framework offers good support for flipping to almost all lectures corresponding to one-fourth of the LA standard curriculum, resulting in 6 videos, 6 to 12 minutes long, per each weekly topic. Then the total structure of videos in one topic mimicking the 3 regular one-hour LA lectures within the university curriculum, but the length of the videos forced the teacher/instructor to be clear and concise, sticking to the kernel of the main ideas and arguments.

The MOOC was carefully designed to pay special attention to various applications of eigenvalues, eigenvectors, and singular values in science and engineering. At the beginning of the MOOC, the eigenvalues are analysed within their central role on geometric transformations, and then they are used to explain the evolution of dynamic systems that include predator-prey models and survival models of species. The topic on singular values will allow the description of a modern application of digital image processing and PCA. From a mathematical point of view, the approach to both these topics includes, in addition to the classic procedure for calculating eigenvalues. and eigenvectors the diagonalization and orthogonal diagonalization of matrices and the spectral theorem. Proofs of main results and theorems can be found in the videos. In the MOOC it is also discussed some numerical results on approximations for eigenvalues and eigenvectors when the matrices are large.

### *2.4. Three years of flipped classroom implementation*

The proposed flipped classroom methodology, described above, was implemented during the last five weeks out of the 14 weeks of three regular first semesters, with improvements between each year.

Being aware of the students' mindset described in section 2.2, the main concerns were: to explain the rationale behind the flipping classroom to students, prepare them for the demanding practice, and supply them with guidelines for orientation during the period of application of the flipped classroom. In the fall semester of 2017/2018, students were weekly reminded (when face-to-face classes took place) of the way that they should watch the videos of the MOOC, practice with quizzes, and put questions in the forum, or bring them to on-campus classes to discuss. Then, during the 2018/2019 a booklet was prepared and made available on the course webpage with guidelines "How to prepare for flipping the classes" to help students to make the most out of the flipped experience. For some students, this was not sufficient, and in the last 2019/2020 experiment, it was also added more information on all communication channels of the course: news in the course webpage about the calendar to follow, weekly individual emails, and post-flipped reports for each class.

During the flipped classroom period, the three weekly lectures (in a total of 12 lectures of one hour long) were changed into an interactive group discussion instead of an expository learning activity centred in the teacher. As it became clear after the first two years of the implementation experiment, the lecture room of "fixed stadium seats" (an auditorium configuration of 250 seats) meant to large classes, was very unsuitable for the implementation of the flipped classroom methodology.

Then, during the 2019/2020 implementation, all efforts were concentrated on moving into the same room (see Fig. 3), where the tutorial/exercise classes took place from the start of the semester. This room, much better suited for blended learning, has a capacity for 40 students (even though, the maximum number of students attending the class was 32). The walls are painted with a special waterproof ink that enables writing on the wall as "on a blackboard". Within this space, the teacher walks amongst the students writing the mathematical formulas on several walls, while the students sit in groups made up of large tables, carrying out the exercises with the freedom to choose using mobile devices, asking a colleague, asking questions to the teacher, looking at the written formulas around, or even doing the exercises directly in some empty spaces on the wall.

Parallel to this, the tutorial classes, i.e. weekly solving exercise sessions, of one hour and a half per week remained very dynamic. Every face-to-face class took place around a table, where eight students had lively discussions about different exercises resolutions, ways to prove mathematical results, after which one designated student per table wrote on the wall the group findings. Then, mediated by the teacher/instructor, all the groups looked around and were invited to comment and reflect on the different mathematical resolutions and methods of proving theorems.

## 3. Results of flipping the class with a MOOC

The validation and monitoring of these three years flipped classroom practices was carried out through quantitative data collection methodologies (final grades of LA, MOOC Técnico platform data, questionnaire surveys, and classes' observations following a protocol) and qualitative data collection methodologies (individual and group interviews, direct and participative observations conducted by the teacher) and detailed results are currently being prepared for a full report. In this preliminary analysis, we'll focus on final grades of LA, success rates on the MOOC, and answers to the questions: "How do you evaluate the flipped classroom methodology?" and "How much the flipped classroom helped you to prepare for the final written test".

### *3.1. Three years of implementation in numbers*

Recall that in each semester, more than 90 CS Engineering first-year students were involved in the flipped classroom experiments. Even if they constitute different groups of students, it's worthwhile to compare their success rates in different years. The final approval rate at LA for students in the 2017/2018 academic year was 74%, in 2018/2019 was 81%, and in the current semester of 2019/2020 is 83%. If we compare these numbers with 67% of approval rate in 2016/2017 the numbers are encouraging us to continue to further develop this kind of learning methodology.

In the online course edition of 2017, from a total of 207 participants, 86 enrollees were first year LA students, from those, 56 completed the online assessments (60% or more of success in graded tests). Meanwhile, in the 2018 MOOC edition of, from 200 participants, 83 enrollees were first-year LA students, and 67 completed the online assessments (60% or more of success in graded tests). And finally, in the 2019 edition, from 149 participants, 88 enrollees were first-year LA students, and 85 completed the online assessments.

From the 45 answers collected through the 2017/2018 final questionnaire directed to the first-year students, we can conclude that regarding the flipped classroom strategy, 70% appreciated or appreciated it a lot (see Fig. 4). Also, 73,3% of the students considered the implementation of this flipping the class model as essential for the preparation of the third on-campus assessment test of LA.

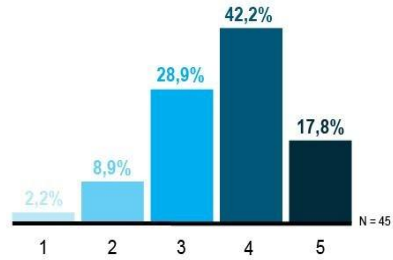

Fig. 4. Satisfaction level of CS students enrolled in 2017/2018 LA course about the flipped classroom methodology (scale: from 1=don't appreciate to 5=appreciated a lot).

Regarding the second-year experience of 2018/2019, from the 48 answers given to the final questionnaire by the first-year students, 56,3% appreciated or appreciated it a lot (see Fig. 5). Besides, 72,9% of the students considered the implementation of this flipping the class model as essential for the preparation of the third on-campus assessment test of LA. Finally, in the third-year experience of 2019/2020, from the 59 answers given to the final questionnaire, 67,7% appreciated or appreciated it a lot (see Fig. 6) and 66,1% claimed that flipped classroom was essential in the preparation for the final assessment test.

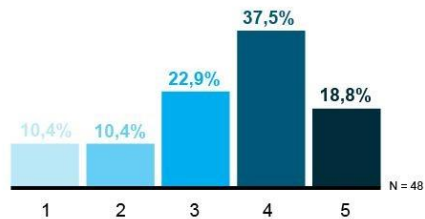

Fig. 5. Satisfaction level of CS students enrolled in 2018/2019 LA course about the flipped classroom methodology (scale: from 1=don't appreciate to 5=appreciated a lot).

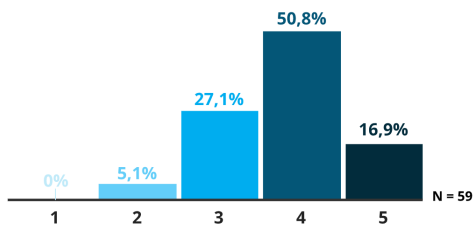

Fig. 6. Satisfaction level of CS students enrolled in 2019/2020 LA course about the flipped classroom methodology (scale: from 1=don't appreciate to 5=appreciated a lot).

*3.2. Discussion of the results*

Based on the completion rates of the online assessment activities, we can identify differences in the attitude between first-year student generations of 2017/2018, 2018/2019, and 2019/2020 semesters. While in the first-generation students were more curious about the learning possibilities of online videos but didn't complete the assessment activities of the MOOC, the second group obtained more success rates in quizzes, and the third group as a whole did complete all online assessments. As a relevant note, we remark that this last group of students also kept very engaged during the face-to-face sessions.

In general, based on the data, we reach the conclusion that the main activity in the MOOC of first-year CS students, all three groups included, consisted of viewing the videos, with the second activity being the discussion about submitting answers to quizzes and other moments of evaluation. Regarding this last design feature (Bates, 2015), we can add that the general feedback received in the open questions, such as general comments, of the final questionnaires, was very helpful. A wide range of assessment activities is assimilated by students, especially in the independent and insecure students, as a great promoter for better content understanding and, at some level, an indirect source of improvement in students' grades.

On the other hand, not having the teacher present, when watching videos at home, to answer immediately to their queries seems to be the biggest perceived disadvantage of inverting the class with virtual learning content. This urgency in the immediate response to questions is the reason why students think that all campus courses can benefit from this flipped classroom strategy, just as a complement to classroom education.

During the three year implementation of the flipped classroom, most of the students who participated in the inverted classes had a positive performance in the written tests. Almost all of them were able to write mathematical statements consistent with careful annotation and respond to conceptual exercises, providing short proofs, if necessary, in addition to algebraic manipulations and sophisticated geometric reasoning. Within each academic year, there was an increase in success rates, both in tests and in questionnaires. At the same time, each group of students recognised more self-confidence and freedom in time management, which helped them in building their learning strategies.

It's worth noting that one student from the first-generation volunteered in the 2018/2019 semester to prepare and conduct one of the flipped classroom sessions. The student himself prepared a topic, with the teacher's prior consent, to initiate the point-to-point discussion, and conducted interactions in the organized classroom. After this action, we understood that the next step could include first-year students preparing guided classes with their colleagues.

## 4. Conclusion and Future work

It is important to note that the work carried out, with the application of the described flipped classroom, in the Linear Algebra course, was recognised as a good practice in 2018 by the IST Observatory of Good Practices.

Even though the practice of the flipped classroom has been, in our case, adopted only in recent years, and with some restrictions, the results are to be considered promising. For first-year students of the university, this challenge provides changes in the way they learn the content from a given subject, in the organisation of their own time, in their own learning strategies and mainly through the development of working group skills. Moreover, the improvements made in the three successive semesters fostered higher students' engagement levels.

It was also clear that the use of the campus room (see Fig. 7) specially adapted for solving group exercises and having group discussions, was crucial for the success of the experiences and the building of a more confident and autonomous group of students. Nevertheless, during COVID-19 confinement pandemic times, this physical room will have to be replaced by dynamic synchronous Zoom meeting rooms, with a mentor/tutor per each group of 8 to 10 students. A teacher or a teaching assistant could play the role of the Zoom's mentor/tutor responsible for the group dynamic.

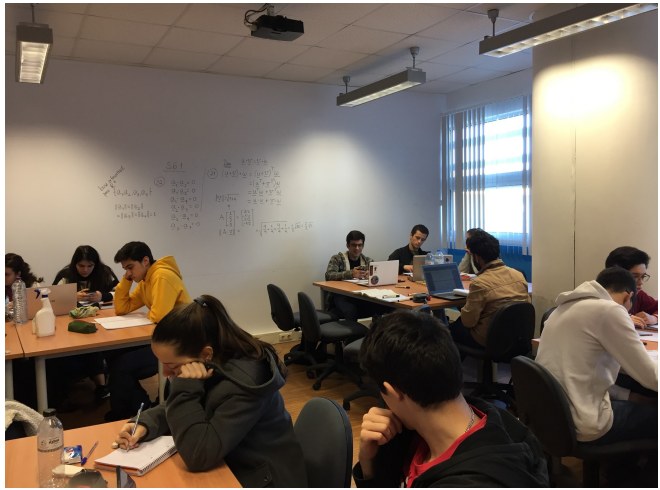

Fig. 7. Students working in groups in the blended-learning campus room. In this room both tutorial and flipped classroom sessions have taken place.

Based on our experience so far, we can say that there is enough information to continue to investin research in this field. We believe that one of the success factors in applying the flipped classroom is having a well-designed MOOC, together with highly motivated instructors and confident students. The MOOC is accessible for everyone at any time to follow and practice, even when one is, due to exceptional circumstances, forced to stay at home, instead of attending face-to-face classes,. Therefore, we have already planned to develop research in the areas of, self-regulation of learning and gamification for future applications of these concepts in other course subjects for which it's possible to run MOOC Técnico online courses.